\documentclass[12pt]{article}
\usepackage[centertags]{amsmath}
\usepackage{amsfonts}
\usepackage{amssymb}
\usepackage{graphicx}
\usepackage{amssymb}
\newtheorem{theorem}{Theorem}[section]

\newtheorem{definition}[theorem]{Definition}

\usepackage[verbose,colorlinks=true,naturalnames=true,linkcolor=blue,]
{hyperref}

\newcounter{rown}

\begin{document}

\title{Quantum Deformations of %% Lorentz and\\ Poincare angebras
Relativistic Symmetries\footnote{Invited talk at the XXII Max Born Symposium "Quantum,
Super and Twistors", September 27-29, 2006 Wroclaw (Poland), in honour of Jerzy
Lukierski.}
\\[10pt]
}
\author{V.N. Tolstoy\footnote{Supported by the grants
RFBR-05-01-01086 and FNRA NT05-241455GIPM.}\\
\\ Institute of Nuclear Physics, Moscow State University, \\
119 992 Moscow, Russia; e-mail: tolstoy@nucl-th.sinp.msu.ru}

%% Quantum Deformations of Relativistic Symmetries
%% (INP, Moscow State University)
%% 10 pages. Invited talk at the XXII Max Born Symposium "Quantum, Super and Twistors",
%% September 27-29, 2006, Wroclaw (Poland), in honour of Jerzy Lukierski.
%% MSC-class: 81R50 (Primary), 17B37 (Secondary)
\date{}

\maketitle
\begin{abstract}
We discussed quantum deformations of $D=4$ Lorentz and Poincar\'{e} algebras. In the case
of Poincar\'{e} algebra it is shown that almost all classical $r$-matrices of
S.~Zakrzewski classification correspond to twisted deformations of Abelian and Jordanian
types. A part of twists corresponding to the $r$-matrices of Zakrzewski classification
are given in explicit form.
\end{abstract}

\section{Introduction}

%% {$22^{th}$ Max Born Symposium\\
%% Wroclaw, Poland, Sept. 27 -- 29,  2006 \\[7pt] %% Date / Occasion
%%Devoted to
%% To my friend Prof. Jurek Lukierski \\ on the occasion of his 70th birthday}
The quantum deformations of relativistic symmetries are described by Hopf-algebraic
deformations of Lorentz and  Poincar\'{e} algebras. Such quantum deformations are
classified by Lorentz and Poincar\'{e} Poisson structures. These Poisson structures given
by classical $r$-matrices were classified already some time ago by S. Zakrzewski in
\cite{Z1} for the Lorentz algebra and in \cite{Z2} for the Poincar\'{e} algebra. In the
case of the Lorentz algebra a complete list of classical $r$-matrices involves the four
independent formulas and the corresponding quantum deformations in different forms were
already discussed in literature (see \cite{ChD, M, BLT1, BLT2, BLT3}). In the case of
Poincar\'{e} algebra the total list of the classical $r$-matrices, which satisfy the
homogeneous classical Yang-Baxter equation, consists of 20 cases which have various
numbers of free parameters. Analysis of these twenty solutions shows that each of them
can be presented as a sum of subordinated $r$-matrices which almost all are of Abelian
and Jordanian types. A part of twists corresponding to the $r$-matrices of Zakrzewski
classification are given in explicit form.

\setcounter{equation}{0}
\section{Preliminaries}

Let $r$ be a classical $r$-matrix of a Lie algebra $\mathfrak{g}$, i.e.
$r\in\,\stackrel{2}\wedge\mathfrak{g}$ and $r$ satisfies to the classical Yang--Baxter
equation (CYBE)
\begin{eqnarray}\label{p1}
%% [[r,r]]\!\!&:=\!\!&
[r^{12},\,r^{13}+\,r^{23}] + [r^{13},\,r^{23}]\!\!&=\!\!&\Omega~,
\end{eqnarray}
where $\Omega$ is $\mathfrak{g}$-invariant element, $\Omega\in(\stackrel{3}\wedge
\mathfrak{g})_{\mathfrak{g}}$. We consider two types of the classical $r$-matrices and
corresponding twists.

Let the classical $r$-matrix $r=r_{A}^{}$ has the form
\begin{eqnarray}\label{p2}
r_{A}^{}\!\!&=\!\!&%%\xi
\sum_{i=1}^{n}x_{i}\wedge y_{i}~,
\end{eqnarray}
where all elements $x_i, y_i$ $(i=1,\ldots,n)$ commute among themselves. Such an
$r$-matrix is called of Abelian type. The corresponding twist is given as follows
\begin{eqnarray}\label{p3}
F_{r_{A}^{}}\!\!&=\!\!&\exp\frac{r_{A}^{}}{2}=
\exp\Bigl(\frac{1}{2}\sum_{i=1}^{n}x_{i}\wedge y_{i}\Bigr)~.
\end{eqnarray}
This twisting two-tensor $F:=F_{r_{A}^{}}$ satisfies the cocycle equation
\begin{equation}\label{p4}
F^{12}(\Delta\otimes{\rm id})(F)\;=\;F^{23}({\rm id}\otimes\Delta)(F)~,
\end{equation}
and the "unital" normalization condition
\begin{equation}\label{p5}
(\epsilon \otimes{\rm id})(F)\;=\;({\rm id}\otimes\epsilon )(F)\;=\;1~.
\end{equation}
The twisting element $F$ defines a deformation of the universal enveloping algebra
$U(\mathfrak{g})$ considered as a Hopf algebra. The new deformed coproduct and antipode
are given as follows
\begin{equation}\label{p6}
\Delta^{(F)}(a)\;=\;F\Delta(a)F^{-1}~,\qquad S^{(F)}(a)=uS(a)u^{-1}
\end{equation}
for any $a\in U(\mathfrak{g})$, where $\Delta(a)$ is a co-product before twisting, and
$u=\sum_i f^{(1)}_{i}S(f^{(2)}_i)$ if $F=\sum_i f^{(1)}_i\otimes f^{(2)}_i$.

Let the classical $r$-matrix $r=r_{J}^{}(\xi)$ has the form\footnote{Here entering the
parameter deformation $\xi$ is a matter of convenience.}
\begin{equation}\label{p7}
r_{J}^{}(\xi)\;=\;\xi\,\Bigl(\sum_{\nu=0}^{n}x_{\nu}\wedge y_{\nu}\Bigr)~,
\end{equation}
where the elements $x_\nu,y_\nu$ $(\nu=0,1,\ldots, n)$ satisfy the relations\footnote{It
is easy to verify that the two-tensor (\ref{p7}) indeed satisfies the homogenous
classical Yang-Baxter equation (\ref{p1}) (with $\Omega=0$), if the elements
$x_\nu,y_\nu$ $(\nu=0,1,\ldots,n)$ are subject to the relations (\ref{p8}).}
\begin{equation}\label{p8}
\begin{array}{rcl}
[x_{0},\,y_{0}]\!\!&=\!\!&y_{0}~,\qquad
[x_{0},\,x_{i}]\;=\;(1-t_{i})x_{i}~,\qquad\;\;[x_{0},\,y_{i}]\;=\;t_{i}y_{i}~,
\\[7pt]
[x_{i},\,y_{j}]\!\!&=\!\!&\delta_{ij}y_{0}~,\quad[x_{i},\,x_{j}]\;=
\;[y_{i},\,y_{j}]\;=\;0~,\quad\;[y_{0},\,x_{j}]\;=\;[y_{0},\,y_{j}]\;=\;0~,
\end{array}
\end{equation}
$(i,j=1,\ldots,n)$, $(t_{i}\in{\mathbb C})$. Such an $r$-matrix is called of Jordanian
type. The corresponding twist is given as follows \cite{T1,T2}
\begin{eqnarray}\label{p9}
F_{r_{J}^{}}\!\!&=\!\!& \exp\Big(\xi\sum\limits_{i=1}^{n}x_{i}\otimes
y_{i}\;e^{-2t_{i}\sigma}\Bigr)\exp(2x_{0}^{}\otimes\sigma)~,
\end{eqnarray}
where $\sigma:=\frac{1}{2}\ln(1+\xi y_{0})$.\footnote{The corresponding twists for Lie
algebras $\mathfrak{sl}(n)$, $\mathfrak{so}(n)$ and $\mathfrak{sp}(2n)$ were firstly
constructed in the papers \cite{KLM, KLO, LSK, AKL}.}

Let $r$ be an arbitrary $r$-matrix of $\mathfrak{g}$. We denote a support of $r$ by
$\mathop{\rm Sup}(r)$\footnote{The support $\mathop{\rm Sup}(r)$ is a subalgebra of
$\mathfrak{g}$ generated by the elements $\{x_i,y_i\}$ if $r=\sum_{i}x_i\wedge y_i$.}.
The following definition is useful.
\begin{definition}
Let $r^{}_1$ and $r^{}_2$ be two arbitrary classical $r$-matrices. We say that $r^{}_2$
is subordinated to $r^{}_1$, $r^{}_1\succ r^{}_2$, if $\delta_{r^{}_1}(\mathop{\rm
Sup}(r^{}_2))=0$, i.e.
\begin{equation}\label{p10}
\delta_{r_{1}^{}}(x)\;:=\;[x\otimes1+1\otimes x,\,r_{1}^{}]\;=\;0~, \quad \forall x\in
\mathop{\rm Sup}(r_2)~.
\end{equation}
\end{definition}
If $r^{}_1\succ r^{}_2$ then $r=r^{}_1+r^{}_2$ is also a classical $r$-matrix (see
\cite{BD}). The subordination enables us to construct a correct sequence of
quantizations. For instance, if the $r$-matrix of Jordanian type (\ref{p7}) is
subordinated to the $r$-matrix of Abelian type (\ref{p2}), $r_{\!A}^{}\succ r_{\!J}^{}$,
then the total twist corresponding to the resulting $r$-matrix $r=r_{\!A}^{}+r_{\!J}^{}$
is given as follows
\begin{eqnarray}\label{p11}
F_{r}\!\!&=\!\!&F_{r_{\!J}^{}}\,F_{r_{\!A}^{}}.
\end{eqnarray}
The further definition is also useful.
\begin{definition}
A twisting two-tensor $F_{r}(\xi)$ of a Hopf algebra, satisfying the conditions
(\ref{p4}) and (\ref{p5}), is called locally $r$-symmetric if the expansion of
$F_{r}(\xi)$ in powers of the parameter deformation $\xi$ has the form
\begin{eqnarray}\label{p12}
F_{r}(\xi)\!\!&=\!\!& 1+c\,r+\mathcal{O}(\xi^2)\ldots
\end{eqnarray}
where $r$ is a classical $r$-matrix, and $c$ is a numerical coefficient, $c\neq0$.
\end{definition}
It is evident that the Abelian twist (\ref{p3}) is globally $r$-symmetric and the twist
of Jordanian type (\ref{p9}) does not satisfy the relation (\ref{p12}), i.e. it is not
locally $r$-symmetric.

\setcounter{equation}{0}
\section{Quantum deformations of Lorentz algebra}

The results of this section in different forms were already discussed in literature (see
\cite{ChD, M, BLT1, BLT2, BLT3}).

The classical canonical basis of the $D=4$ Lorentz algebra, $\mathfrak{o}(3,1)$, can be
described by anti-Hermitian six generators ($h$, $e_{\pm}$, $h'$, $e'_{\pm}$) satisfying
the following non-vanishing commutation relations\footnote{Since the real Lie algebra
$\mathfrak{o}(3,1)$ is standard realification of the complex Lie
$\mathfrak{sl}(2,\mathbb{C})$ these relations are easy obtained from the defining
relations for $\mathfrak{sl}(2,\mathbb{C})$, i.e. from (\ref{L1}).}:
\begin{eqnarray}
&[h,\,e_{\pm}^{}]\;=\;\pm e_{\pm}^{}\,,\qquad [e_{+}^{},\,e_{-}^{}]\;=\;2h~,
\label{L1}
\\[8pt]
&[h,\,e'_{\pm}]\;=\;\pm e'_{\pm}~,\qquad [h',\,e_{\pm}]\;=\;\pm e'_{\pm}~,\qquad
[e_{\pm}^{},\,e'_{\mp}]\;=\;\pm2h'~ ,
\label{L2}
\\[8pt]
& [h',\,e'_{\pm}]\;=\;\mp e_{\pm}^{}~,\qquad [e'_{+},\,e'_{-}]\;=\;-2h~,
\label{L3}
\end{eqnarray}
and moreover
\begin{equation}\label{L4}
x^*\;=\;-x\qquad (\forall\;x\;\in\; \mathfrak{o}(3,1))~.\textbf{}
\end{equation}
A complete list of classical $r$-matrices which describe all Poisson structures and
generate quantum deformations for $\mathfrak{o}(3,1)$ involve the four independent
formulas \cite{Z1}:
\begin{eqnarray}\label{L5}
r_{1}^{}\!\!&=\!\!&\alpha\,e_{+}\wedge h~,
\\[5pt]\label{L6}
r_{2}^{}\!\!&=\!\!&\alpha\,(e_{+}\wedge h-e'_{+}\wedge h')+2\beta\, e'_{+}\wedge e_{+}~,
\\[5pt]\label{L7}
r_{3}^{}\!\!&=\!\!&\alpha\,(e'_{+}\wedge e_{-}+ e_{+}\wedge e'_{-})\,+\,
\beta\,(e_{+}\wedge e_{-}\,-\,e'_{+}\wedge e'_{-})-2\gamma\,h\wedge h'~,
\\[5pt]\label{L8}
r_{4}^{}\!\!&=\!\!&\alpha\bigl(e'_{+}\wedge e_{-}+ e_{+}\wedge e'_{-}-2h\wedge
h'\bigr)\pm e_{+}\wedge e'_{+}~.
\end{eqnarray}
If the universal $R$-matrices of the quantum deformations corresponding to the classical
$r$-matrices (\ref{L5})--(\ref{L8}) are unitary then these $r$-matrices are
anti-Hermitian, i.e.
\begin{equation}\label{L9}
r^*_{j}\;=\;-r_{j}\qquad (j=1,2,3,4)~.
\end{equation}
Therefore the $*$-operation (\ref{L4}) should be lifted to the tensor product
$\mathfrak{o}(3,1)\otimes \mathfrak{o}(3,1)$. There are two variants of this lifting:
{\it direct} and {\it flipped}, namely,
\begin{eqnarray}\label{L10}
(x\otimes y)^*\!\!&=\!\!&x^*\otimes y^*\qquad({\rm *-direct})~,
\\[5pt]\label{L11}
(x\otimes y)^*\!\!&=\!\!&y^*\otimes x^*\qquad({\rm *-flipped})~.
\end{eqnarray}
We see that if the "direct" lifting of the $*$-operation (\ref{L4}) is used then all
parameters in (\ref{L5})--(\ref{L8}) are pure imaginary. In the case of the "flipped"
lifting (\ref{L11}) all parameters in (\ref{L5})--(\ref{L8}) are real.

The first two $r$-matrices (\ref{L5}) and (\ref{L6}) satisfy the homogeneous CYBE and
they are of Jordanian type. If we assume (\ref{L10}), the corresponding quantum
deformations were described detailed in the paper \cite{BLT2} and they are entire defined
by the twist of Jordanian type:
\begin{equation}\label{L12}
F_{r_1^{}}^{}\,=\,\exp{(h\otimes\sigma})~,\quad\sigma\,=\,\frac{1}{2}\ln(1+\alpha e_{+})
\end{equation}
for the $r$-matrix (\ref{L5}), and
\begin{equation}\label{L13}
F_{r_2^{}}^{}\,=\,\exp{\Bigl(\frac{\imath\beta}{\alpha^2}\;\sigma\wedge\varphi\Bigr)}\,
\exp{(h\otimes\sigma-h'\otimes\varphi)}~,
\end{equation}
\begin{equation}\label{L14}
\sigma \,=\,\frac{1}{2}\ln\left[(1+\alpha e_+)^2\!+(\alpha e'_+)^2
\right],\quad\varphi\,=\,\arctan{\frac{\alpha e'_+}{1+\alpha e_+}}
\end{equation}
for the $r$-matrix (\ref{L6}). It should be recalled that the twists (\ref{L12}) and
(\ref{L13}) are not locally $r$-symmetric. A locally $r$-symmetric twist for the
$r$-matrix (\ref{L5}) was obtained in \cite{Oh} and it has the following complicated
formula:
\begin{equation}\label{L15}
F_{r_1^{}}'=\,\exp{\Bigl(\frac{1}{2}\Delta(h)-\frac{1}{2}\Bigl(h \frac{\sinh\alpha
e_{+}}{\alpha e_{+}}\otimes e^{-\alpha e_{+}}\!+ e^{\alpha e_{+}}\!\otimes
h\frac{\sinh\alpha e_{+}}{\alpha e_{+}}\Bigr)
\frac{\alpha\Delta(e_{+})}{\sinh\alpha\Delta(e_{+})}\,\Bigr)},
\end{equation}
where $\Delta$ is a primitive coproduct.

The last two $r$-matrices (\ref{L7}) and (\ref{L8}) satisfy the non-homogeneous
(modified) CYBE and they can be easily obtained from the solutions of the complex algebra
$\mathfrak{o}(4,\mathbb{C})\simeq \mathfrak{sl}
(2,\mathbb{C})\oplus\mathfrak{sl}(2,\mathbb{C})$ which describes the  complexification of
$\mathfrak{o}(3,1)$. Indeed, let us introduce the complex basis of Lorentz algebra
$(\mathfrak{o}(3,1)\simeq\mathfrak{sl}(2;\mathbb{C})\oplus\mathfrak
{\overline{sl}}(2,\mathbb{C}))$ described by two commuting sets of complex generators:
\begin{eqnarray}\label{L16}
H_1\!\!&=\!\!& \frac{1}{2}\,(h+\imath h')~,\qquad
E_{1\pm}\;=\;\frac{1}{{2}}\,(e_{\pm}^{}+ \imath e'_{\pm})~,
\\[5pt]\label{L17}
H_2\!\!&=\!\!& \frac{1}{{2}}\, (h-\imath h')~,\qquad
E_{2\pm}\;=\;\frac{1}{{2}}\,(e_{\pm}^{}-\imath e'_{\pm})~,
\end{eqnarray}
which satisfy the relations (compare with (\ref{L1}))
\begin{equation}\label{L18}
[H_k,\,E_{k\pm}]\;=\;\pm E_{k\pm}~,\qquad [E_{k+},\,E_{k-}]\;=\;2 H_k\qquad(k=1,2)~.
\end{equation}
The $*$-operation describing the real structure acts on the generators $H_k$, and
$E_{k\pm}$ ($k=1,2$) as follows
\begin{equation}\label{L19}
H_1^*\;=\;-H_2^{}~,\qquad E_{1\pm}^*\;=\;- E_{2\pm}^{}~,\qquad H_2^*\;=\;-H_1^{}~,\qquad
E_{2\pm}^*\;=\;- E_{1\pm}^{}~.
\end{equation}
The classical $r$-matrix $r_3$, (\ref{L7}), and $r_4$, (\ref{L8}), in terms of the
complex basis (\ref{L16}), (\ref{L17}) take the form
\begin{equation}\label{L20}
\begin{array}{rcl}
r_3^{}\!\!&=\!\!&r'_1+r''_1~,
\\[7pt]
r'_3\!\!&:=\!\!&2(\beta+\imath\alpha)E_{1+}\wedge
E_{1-}+2(\beta-\imath\alpha)E_{2+}\wedge E_{2-}~,
\\[7pt]
r''_3\!\!&:=\!\!&4\imath\gamma\,H_{2}\wedge H_{1}~,
\end{array}
\end{equation}
and
\begin{equation}\label{L21}
\begin{array}{rcl}
r_4^{}\!\!&=\!\!&r'_4+r''_4~,\quad
\\[7pt]
r'_4\!\!&:=\!\!&2\imath\alpha(E_{1+}\wedge E_{1-}-E_{2+}\wedge E_{2-}-2H_{1}\wedge
H_{2})~,
\\[7pt]
r''_4\!\!&:=\!\!&4\imath\lambda\,E_{1+}\wedge E_{2+} %% \quad(\lambda\in\mathbb{R})~.
\end{array}
\end{equation}

For the sake of convenience we introduce parameter\footnote{We can reduce this parameter
$\lambda$ to $\pm\frac{1}{2}$ by automorphism of $\mathfrak{o}(4,\mathbb{C})$.}$\lambda$
in $r''_{4}$. It should be noted that $r'_{3}$, $r''_{3}$ and $r'_{4}$, $r''_{4}$ are
themselves classical $r$-matrices. We see that the $r$-matrix $r'_{3}$ is simply a sum of
two standard $r$-matrices of $\mathfrak{sl}(2;\mathbb{C})$, satisfying the anti-Hermitian
condition $r^*=-r$. Analogously, it is not hard to see that the $r$-matrix $r_4$
corresponds to a Belavin-Drinfeld triple \cite{BD} for the Lie algebra
$\mathfrak{sl}(2;\mathbb{C})\oplus \mathfrak{\overline{sl}}(2,\mathbb{C}))$. Indeed,
applying the Cartan automorphism $E_{2\pm}\rightarrow E_{2\mp}$, $H_{2}\rightarrow-H_{2}$
we see that this is really correct (see also \cite{IO2001}).

We firstly describe quantum deformation corresponding to the classical $r$-matrix $r_{3}$
(\ref{L20}). Since the $r$-matrix $r''_{3}$ is Abelian and it is subordinated to $r'_{3}$
therefore the algebra $\mathfrak{o}(3,1)$ is firstly quantized in the direction $r'_{3}$
and then an Abelian twist corresponding to the $r$-matrix $r''_{3}$ is applied. We
introduce the complex notations $z_{\pm}:=\beta\pm\imath\alpha$. It should be noted that
$z_{-}^{}=z_{+}^*$ if the parameters $\alpha$ and $\beta$ are real, and
$z_{-}^{}=-z_{+}^*$ if the parameters $\alpha$ and $\beta$ are pure imaginary. From
structure of the classical $r$-matrix $r_{3}'$ it follows that a quantum deformation
$U_{r'_{1}} (\mathfrak{o}(3,1))$ is a combination of two $q$-analogs of
$U(\mathfrak{sl}(2;\mathbb{C}))$ with the parameter $q_{z_{+}^{}}$ and $q_{z_{-}^{}}$,
where $q_{z_{\pm}^{}}:=\exp{z_{\pm}^{}}$. Thus $U_{r'_{3}}(\mathfrak{o}(3,1))\cong
U_{q_{z_{+}^{}}^{}}(\mathfrak{sl}(2;\mathbb{C}))\otimes U_{q_{z_{-}^{}}^{}}
(\overline{\mathfrak{sl}}(2;\mathbb{C}))$ and the standard generators $q_{z_{+}^{}}^{\pm
H_{1}}$, $E_{1\pm}$ and $q_{z_{-}^{}}^{\pm H_{2}}$, $E_{2\pm}$ satisfy the following
non-vanishing defining relations
\begin{eqnarray}\label{L22}
q_{z_{+}}^{H_1}E_{1\pm}\!\!&=\!\!&q_{z_{+}}^{\pm1}E_{1\pm}\,q_{z_{+}}^{H_1}, \qquad
[E_{1+},\,E_{1-}]\;=\;\frac{q_{z_{+}}^{2H_1}-q_{z_{+}}^{-2H_1}}
{q_{z_{+}}^{}-q_{z_{+}}^{-1}},
\\[7pt]\label{L23}
q_{z_{-}}^{H_2}E_{2\pm}\!\!&=\!\!&q_{z_{-}}^{\pm1}E_{2\pm}\,q_{z_{-}}^{H_2}, \qquad
[E_{2+},\,E_{2-}]\;=\;\frac{q_{z_{-}}^{2H_2}-q_{z_{-}}^{-2H_2}}
{q_{z_{-}}^{}-q_{z_{-}}^{-1}}~.
\end{eqnarray}
In this case the co-product $\Delta_{r'_{1}}$ and antipode $S_{r'_{1}}$ for the
generators $q_{z_{+}^{}}^{\pm H_{1}}$, $E_{1\pm}$ and $q_{z_{-}^{}}^{\pm H_{2}}$,
$E_{2\pm}$ can be given by the formulas:
\begin{eqnarray}\label{L24}
\Delta_{r'_{1}}^{}(q_{z_{+}}^{\pm H_{1}})\!\!&=\!\!&q_{z_{+}}^{\pm H_{1}}\otimes
q_{z_{+}}^{\pm H_{1}}, \qquad\Delta_{r'_{1}}^{}(E_{1\pm})\;=\;E_{1\pm}\otimes
q_{z_{+}}^{H_{1}}+q_{z_{+}}^{-H_{1}}\otimes E_{1\pm}~,
\\[7pt]\label{L25}
\Delta_{r'_{1}}^{}(q_{z_{-}}^{\pm H_{2}})\!\!&=\!\!&q_{z_{-}}^{\pm H_{2}}\otimes
q_{z_{-}}^{\pm H_{2}},\qquad\Delta_{r'_{1}}^{}(E_{2\pm})\;=\;E_{2\pm}\otimes
q_{z_{-}}^{H_{2}}+q_{z_{-}}^{-H_{2}}\otimes E_{2\pm}~,
\end{eqnarray}
\begin{eqnarray}\label{L26}
S_{r'_{1}}^{}(q_{z_{+}}^{\pm H_{1}})\!\!&=\!\!&q_{z_{+}}^{\mp H_{1}},\qquad
S_{r'_{1}}^{}(E_{1\pm})\;=\;-q_{z_{+}}^{\pm1}E_{1\pm}~,
\\[7pt]\label{L27}
S_{r'_{1}}^{}(q_{z_{-}^{}}^{\pm H_{2}})\!\!&=\!\!&q_{z_{-}}^{\mp H_{2}},\qquad
S_{r'_{1}}^{}(E_{2\pm})\;=\;-q_{z_{-}}^{\,\pm1}E_{2\pm}~.
\end{eqnarray}
The $*$-involution describing the real structure on the generators (\ref{L8}) can be
adapted to the quantum generators as follows
\begin{equation}\label{L28}
(q_{z_{+}}^{\pm H_{1}})^*\;=\;q_{z_{+}^{*}}^{\mp H_{2}},\quad\;E_{1\pm}^*\;=\;-
E_{2\pm}^{}~,\quad\;(q_{z_{-}}^{\pm H_{2}})^*\;=\;q_{z_{-}^{*}}^{\mp H_{1}},
\quad\;E_{2\pm}^*\;=\;- E_{1\pm}^{}~,
\end{equation}
and there exit two $*$-liftings: {\it direct} and {\it flipped}, namely,
\begin{eqnarray}\label{L28}
(a\otimes b)^*\!\!&=\!\!&a^*\otimes b^*\qquad({\rm*-direct})~,
\\[5pt]\label{L30}
(a\otimes b)^*\!\!&=\!\!&b^*\otimes a^*\qquad({\rm*-flipped})
\end{eqnarray}
for any $a\otimes b\in U_{r'_{3}}(\mathfrak{o}(3,1))\otimes U_{r'_{3}}
(\mathfrak{o}(3,1))$, where the $*$-direct involution corresponds to the case of the pure
imaginary parameters $\alpha,\,\beta$ and the $*$-flipped involution corresponds to the
case of the real deformation parameters $\alpha,\,\beta$.  It should be stressed that the
Hopf structure on $U_{r'_{3}}(\mathfrak{o}(3,1))$ satisfy the consistency conditions
under the $*$-involution
\begin{equation}\label{L31}
\Delta_{r'_{3}}(a^*)\;=\;(\Delta_{r'_{3}}(a))^*,\quad\;\; S_{r'_{3}}((S_{r'_{3}}
(a^*))^{*})\;=\;a\quad(\forall x\in U_{r'_{3}}(\mathfrak{o}(3,1))~.
\end{equation}

Now we consider deformation of the quantum algebra $U_{r'_{3}}(\mathfrak{o}(3,1))$
(secondary quantization of $U(\mathfrak{o}(3,1))$) corresponding to the additional
$r$-matrix $r''_{3}$, (\ref{L20}). Since the generators $H_{1}$ and $H_{2}$ have the
trivial coproduct
\begin{eqnarray}\label{L32}
\Delta_{r'_{3}}(H_{k})\!\!&=\!\!&H_{k}\otimes 1+1\otimes H_{k}\quad(k=1,2)~,
\end{eqnarray}
therefore the unitary two-tensor
\begin{eqnarray}\label{L33}
F_{r_{3}''}^{}\!\!:=\!\!&q_{\imath\gamma}^{H_{1}\wedge H_{2}}\qquad
(F_{r_{3}''}^*\;=\;F_{r_{1}''}^{-1})
\end{eqnarray}
satisfies the cocycle condition (\ref{p4}) and the "unital" normalization condition
(\ref{p5}). Thus the complete deformation corresponding to the $r$-matrix $r_{3}^{}$ is
the twisted deformation of $U_{r'_{3}}(\mathfrak{o}(3,1))$, i.e. the resulting coproduct
$\Delta_{r_{3}^{}}$ is given as follows
\begin{eqnarray}\label{L34}
\Delta_{r_{3}^{}}^{}(x)\!\!&=\!\!&F_{r_{1}''}^{}\Delta_{r_{1}'}^{}(x)
F_{r_{3}''}^{-1}\quad(\forall x\in U_{r'_{1}}(\mathfrak{o}(3,1))~.
\end{eqnarray}
and in this case the resulting antipode $S_{r_{3}^{}}^{}$ does not change,
$S_{r_{3}^{}}^{}=S_{r'_{3}}^{}$. Applying the twisting two-tensor (\ref{L33}) to the
formulas (\ref{L24}) and (\ref{L25}) we obtain
\begin{eqnarray}\label{L35}
\Delta_{r_{3}^{}}(q_{z_{+}}^{\pm H_{1}})\!\!&=\!\!&q_{z_{+}}^{\pm H_{1}}\otimes
q_{z_{+}}^{\pm H_{1}},\quad\Delta_{r'_{1}}(q_{z_{-}}^{\pm H_{2}})\;=\;q_{z_{-}}^{\pm
H_{2}}\otimes q_{z_{-}}^{\pm H_{2}},
\\[7pt]\label{L36}
\Delta_{r_{3}^{}}(E_{1+})\!\!&=\!\!&E_{1+}\otimes q_{z_{+}}^{H_{1}}
q_{\imath\gamma}^{H_{2}}+q_{z_{+}}^{-H_{1}} q_{\imath\gamma}^{-H_{2}}\otimes E_{1+}~,
\\[7pt]\label{L37}
\Delta_{r_{3}^{}}(E_{1-})\!\!&=\!\!&E_{1-}\otimes q_{z_{+}}^{H_{1}}
q_{\imath\gamma}^{-H_{2}}+q_{z_{+}}^{-H_{1}}q_{\imath\gamma}^{H_{2}}\otimes E_{1-}~,
\\[7pt]\label{L38}
\Delta_{r_{3}}(E_{2+})\!\!&=\!\!&E_{2+}\otimes q_{z_{-}}^{H_{2}}
q_{\imath\gamma}^{-H_{1}}+q_{z_{-}}^{-H_{2}} q_{\imath\gamma}^{H_{1}}\otimes E_{2+}~,
\\[7pt]\label{L39}
\Delta_{r_{3}}(E_{2-})\!\!&=\!\!&E_{2-}\otimes q_{z_{-}}^{H_{2}}
q_{\imath\gamma}^{H_{1}}+q_{z_{-}}^{-H_{2}}q_{\imath\gamma}^{-H_{1}}\otimes E_{2-}~.
\end{eqnarray}

Next, we describe quantum deformation corresponding to the classical $r$-matrix $r_{4}$
(\ref{L21}). Since the $r$-matrix $r_{4}'(\alpha):=r_4'$ is a particular case of
$r_{3}^{}(\alpha,\beta,\gamma):=r_{3}^{}$, namely
$r_{4}'(\alpha)=r_{3}^{}(\alpha,\beta=0,\gamma=\alpha)$, therefore a quantum deformation
corresponding to the $r$-matrix $r_4'$ is obtained from the previous case by setting
$\beta=0,\gamma=\alpha$, and we have the following formulas for the coproducts
$\Delta_{r_{4}'}$:
\begin{eqnarray}\label{L40}
\Delta_{r_{4}'}(q_{\xi}^{\pm H_{k}})\!\!&=\!\!&q_{\xi}^{\pm H_{k}}\otimes q_{\xi}^{\pm
H_{k}} \qquad (k=1,2)~,
%% \quad\Delta_{r'_{1}}(q_{z_{-}}^{\pm H_{2}})\;=\;q_{z_{-}}^{\pm H_{2}}\otimes
%% q_{z_{-}}^{\pm H_{2}},
\\[7pt]\label{L41}
\Delta_{r_{4}'}(E_{1+})\!\!&=\!\!&E_{1+}\otimes q_{\xi}^{H_{1}+H_{2}}+
q_{\xi}^{-H_{1}-H_{2}}\otimes E_{1+}~,
\\[7pt]\label{L42}
\Delta_{r_{4}'}(E_{1-})\!\!&=\!\!&E_{1-}\otimes q_{\xi}^{H_{1}-H_{2}}+
q_{\xi}^{-H_{1}+H_{2}}\otimes E_{1-}~,
\\[7pt]\label{L43}
\Delta_{r_{4}'}(E_{2+})\!\!&=\!\!&E_{2+}\otimes q_{\xi}^{-H_{1}-H_{2}}+
q_{\xi}^{H_{1}+H_{2}}\otimes E_{2+}~,
\\[7pt]\label{L44}
\Delta_{r_{4}'}(E_{2-})\!\!&=\!\!&E_{2-}\otimes q_{\xi}^{H_{1}-H_{2}}+
q_{\xi}^{-H_{1}+H_{2}}\otimes E_{2-}~,
\end{eqnarray}
where we set $\xi:=\imath\alpha$.

Consider the two-tensor
\begin{eqnarray}\label{L45}
F_{r_{4}''}\!\!:=\!\!&\exp_{q^{2}}^{}\big(\lambda E_{1+}q_{\xi}^{H_{1}+H_{2}}\otimes
E_{2+}q_{\xi}^{H_{1}+H_{2}}\big)~.
\end{eqnarray}
Using properties of $q$-exponentials (see \cite{KT1}) is not hard to verify that
$F_{r_{4}''}$ satisfies the cocycle equation (\ref{p4}). Thus the quantization
corresponding to the $r$-matrix $r_4$ is the twisted $q$-deformation
$U_{r_{4}'}(\mathfrak{o}(3,1))$. Explicit formulas of the co-products
$\Delta_{r_{4}^{}}^{}(\cdot)=F_{r_{4}''}^{}\Delta_{r_{4}'}^{}(\cdot) F_{r_{4}''}^{-1}$
and antipodes $S_{r_4}(\cdot)$ in the complex and real Cartan-Weyl bases of
$U_{r_{4}'}(\mathfrak{o}(3,1))$ will be presented in the outgoing paper \cite{BLT3}.

 \setcounter{equation}{0}
\section{Quantum deformations of Poincare algebra}

The Poincar\'{e} algebra ${\mathcal{P}}(3,1)$ of the 4-dimensional space-time is
generated by 10 elements: the six-dimensinal Lorentz algebra
$\mathfrak{o}(3,1)$ with the generators  $M_i$, $N_i$ ($i=1,2,3$):%%\\[3pt]
\begin{equation}\label{P1}
[M_i,\,M_j ]\;=\;\imath\epsilon_{ijk}\,M_k,\;\;
[M_i,\,N_j]\;=\;\imath\epsilon_{ijk}\,N_k,\;\;[N_i,\,N_j]\;=\;-\imath
\epsilon_{ijk}\,M_k,
\end{equation}
and the four-momenta $P_j$, $P_0$ $(j=1,2,3)$ with the standard commutation relations:
\begin{eqnarray}\label{P2}
[M_j,\,P_k]\!\!&=\!\!&\imath\epsilon_{jkl}\,P_l~,\qquad [M_j,\,P_0]\;=\;0~,
\\[5pt]\label{P3}
[N_j,\,P_k]\!\!&=\!\!&-\imath\delta_{jk}\,P_0~,\quad\; [N_j,\,P_0]\;=\;-\imath P_j^{}~.
\end{eqnarray}
The physical generators of the Lorentz algebra, $M_i$, $N_i$ ($i = 1,2,3$), are related
with the canonical basis $h,h',e_{\pm},e'_{\pm}$ as follows
\begin{eqnarray}\label{P4}
h\!\!&=\!\!&\imath N_3~,\qquad\quad e_{\pm}\;=\;\imath (N_1\pm\,M_2),
\\[5pt]\label{P5}
h'\!\!&=\!\!&-\imath M_3~,\qquad e'_{\pm}\;=\;\imath (\pm N_2-M_1).
\end{eqnarray}
The subalgebra generated by the four-momenta $P_j$, $P_0$ $(j=1,2,3)$ will be denoted by
$\mathbf{P}$ and we also set $P_{\pm}:=P_{0}\pm P_{3}$.

S.~Zakrzewski has shown in \cite{Z2} that each classical $r$-matrix,
$r\in\mathcal{P}(3,1) \wedge\mathcal{P}(3,1)$,  has a decomposition
\begin{equation}\label{P6}
r=a+b+c~,
\end{equation}
where $a\in\mathbf{P}\wedge\mathbf{P}$, $b\in\mathfrak{o}(3,1)\wedge\mathbf{P}$,
$c\in\mathfrak{o}(3,1)\wedge\mathfrak{o}(3,1)$ satisfy the following relations
\begin{eqnarray}\label{P7}
[[c,c]]\!\!&=\!\!&0~,
\\[3pt]\label{P8}
[[b,c]]\!\!&=\!\!&0~,
\\[3pt]\label{P9}
2[[a,c]]+[[b,b]]\!\!&=\!\!&t\Omega\quad (t\in \mathbb{R})~,
\\[3pt]\label{P10}
[[a,b]]\!\!&=\!\!&0~.
\end{eqnarray}
Here $[[\cdot,\cdot]]$ means the Schouten bracket. Moreover a total list of the classical
$r$-matrices for the case $c\neq0$ and also for the case $c=0$, $t=0$ was
found.\footnote{Classification of the $r$-matrices for the case $c=0$, $t\neq0$ is an
open problem up to now.} It was shown that there are fifteen solutions for the case
$c=0$, $t=0$, and six solutions for the case $c\neq0$ where there is only one solution
for $t\neq0$. Thus Zakrzewski found twenty $r$-matrices which satisfy the homogeneous
classical Yang-Baxter equation ($t=0$ in (\ref{P9})). Analysis of these twenty solutions
shows that each of them can be presented as a sum of subordinated $r$-matrices which
almost all  are of Abelian and Jordanian types. Therefore these $r$-matrices correspond
to twisted deformations of the Poincar\'{e} algebra ${\mathcal{P}}(3,1)$. We present here
$r$-matrices only for the case $c\neq0,\;t=0$:
%%\footnote{The rest solutions which correspond to the case $c=0,\;t=0$ are simpler.}
\begin{eqnarray}\label{P11}
r_1\!\!&=\!\!&\gamma h'\wedge h+\alpha(P_{+}\wedge P_{-}-P_{1}\wedge P_{2})~,
\\[5pt]\label{P12}
r_2\!\!&=\!\!&\gamma e'_{+}\wedge e_{+}+\beta_{1}(e_{+}\wedge P_{1}-e'_{+}\wedge P_{2}+
h\wedge P_{+})+\beta_{2}h'\wedge P_{+}~,
\\[5pt]\label{P13}
r_3\!\!&=\!\!&\gamma e'_{+}\wedge e_{+}+\beta(e_{+}\wedge P_{1}-e'_{+}\wedge P_{2}
+h\wedge P_{+})+\alpha P_{1}\wedge P_{+}~,
\\[5pt]\label{P14}
r_4\!\!&=\!\!&\gamma(e'_{+}\wedge e_{+}+e_{+}\wedge P_{1}\!+e'_{+}\wedge P_{2}\!-
P_{1}\wedge P_{2})+P_{+}\wedge(\alpha_{1}P_{1}\!+\alpha_{2}P_{2})~,
\\[5pt]\label{P15}
r_5\!\!&=\!\!&\gamma_{1}(h\wedge e_{+}-h'\wedge e'_{+})+ \gamma_{2}e_{+}\wedge e'_{+}~.
\end{eqnarray}

The first $r$-matrix $r_1$ is a sum of two subordinated Abelian $r$-matrices
\begin{eqnarray}\label{P16}
\begin{array}{rcl}
r_1\!\!&:=\!\!&r_1'+r_1''~,\quad r_1'\,\succ\,r_1''~,
\\[5pt]
r_1'\!\!&=\!\!&\alpha(P_{+}\wedge P_{-}-P_{1}\wedge P_{2})~,\qquad r_1''\;:=\;\gamma
h'\wedge h~.
\end{array}
\end{eqnarray}
Therefore the total twist defining quantization in the direction to this $r$-matrix is
the ordered product of  two the Abelian twits
\begin{eqnarray}\label{P17}
F_{r_1^{}}\!\!&=\!\!&F_{r_1''}F_{r_1'}\;=\;\exp\bigl(\gamma h'\wedge h\bigr)\,
\exp\bigl(\alpha(P_{+}\wedge P_{-}-P_{1}\wedge P_{2}) \bigr)~.
\end{eqnarray}

The second  $r$-matrix $r_2$ is a sum of three subordinated $r$-matrices where two of
them are of Abelian type and one is of Jordanian type
\begin{eqnarray}\label{P18}
\begin{array}{rcl}
r_2\!\!&:=\!\!&r_{3}'+r_{2}''+r_{2}'''~,\quad r_{2}'\,\succ\, r_{2}''\,\succ\, r_{2}'''~,
\\[5pt]
r_{2}'\!\!&:=\!\!&\beta_{1}(e_{+}\wedge P_{1}-e'_{+}\wedge P_{2}+h\wedge P_{+})~,
\\[5pt]
r_{2}''\!\!&:=\!\!&\gamma e'_{+}\wedge e_{+}~,\quad r_{2}'''\;:=\;\beta_{2} h'\wedge
P_{+}~.
\end{array}
\end{eqnarray}
Corresponding twist is given by the following formulas
\begin{eqnarray}\label{P19}
F_{r_2^{}}\!\!&=\!\!&F_{r_{2}'''}F_{r_{2}''}F_{r_{2}'}~,
\end{eqnarray}
where
\begin{eqnarray}\label{P20}
\begin{array}{rcl}
F_{r_{2}^{'}}\!\!&=\!\!& \exp\bigr(\beta_{1}(e_{+}^{}\otimes P_{1}- e_{+}'\otimes
P_{2})\bigr)\exp(2h\otimes\sigma_{+})~,
\\[10pt]  %%\label{P20}
F_{r_{2}^{''}}\!\!&=\!\!&\exp(\gamma e_{+}'\wedge e_{+})~,\quad
F_{r_{2}'''}\;=\;\exp(\beta_2h'\wedge \sigma_{+})~.
%% \\[10pt]
%% \tilde{h}'\!\!&=\!\!&\omega h'\omega^{-1}~,\qquad \displaystyle\omega\,=\,\exp\Bigl(
%% \frac{-\beta_{1}\sigma_{+}}{\exp2\sigma_{+}-1}(e_{+}^{}P_{1}- e_{+}'P_{2})\Bigr).
\end{array}
\end{eqnarray}
Here and below we set $\sigma_+:=\frac{1}{2}\ln(1+\beta_{1}P_{+})$.

The third $r$-matrix $r_3$ is a sum of two subordinated $r$-matrices where one is of
Abelian type and another is a more complicated $r$-matrix which we call mixed
Jordanian-Abelian type
\begin{eqnarray}\label{P21}
\begin{array}{rcl}
r_3\!\!&:=\!\!&r_{3}'+r_{3}''~,\quad r_{3}'\;\succ\, r_{3}''~,
\\[5pt]
r_{3}'\!\!&:=\!\!&\beta_{1}(e_{+}\wedge P_{1}-e'_{+}\wedge P_{2}+h\wedge P_{+})+\alpha
P_{1}\wedge P_{+}~,
\\[5pt]
r_{3}''\!\!&:=\!\!&\gamma e'_{+}\wedge e_{+}~.
\end{array}
\end{eqnarray}
Corresponding twist is given by the following formulas
\begin{eqnarray}\label{P22}
F_{r_3^{}}\!\!&=\!\!&F_{r_{3}''}F_{r_{3}'}~,
\end{eqnarray}
where
\begin{eqnarray}\label{P23}
\begin{array}{rcl}
F_{r_{3}^{'}}\!\!&=\!\!& \exp\bigr(\beta_{1}(e_{+}^{}\otimes P_{1}- e_{+}'\otimes
P_{2})\bigr)\exp(\alpha P_{1}\wedge\sigma_{+})\exp(2h\otimes\sigma_{+})~,
\\[10pt]  %
F_{r_{3}^{''}}\!\!&=\!\!&\exp(\gamma e_{+}'\wedge e_{+})~.
\end{array}
\end{eqnarray}

The fourth $r$-matrix $r_4$ is a sum of two subordinated $r$-matrices of Abelian type
\begin{eqnarray}\label{P24}
\begin{array}{rcl}
r_4\!\!&:=\!\!&r_{4}'+r_{4}''~,\quad r_{4}'\;\succ\, r_{4}''~,
\\[5pt]
r_{4}'\!\!&:=\!\!&P_{+}\wedge(\alpha_{1}P_{1}+\alpha_{+}P_{2})~,
\\[5pt]
r_{4}''\!\!&:=\!\!&\gamma(e'_{+}-P_{1})\wedge(e_{+}^{}+P_{2})~.
\end{array}
\end{eqnarray}
Corresponding twist is given by the following formulas
\begin{eqnarray}\label{P22}
F_{r_4^{}}\!\!&=\!\!&F_{r_{4}''}F_{r_{4}'}~,
\end{eqnarray}
where
\begin{eqnarray}\label{P23}
\begin{array}{rcl}
F_{r_{4}^{'}}\!\!&=\!\!& \exp\bigr((P_{+}^{}\otimes (\alpha_{1}P_{1}+
\alpha_{2}P_{2})\bigr)~,
\\[5pt]
F_{r_{4}^{''}}\!\!&=\!\!& \exp\bigr(\gamma(e'_{+}-P_{1})\wedge(e_{+}^{}+P_{2})\bigr)~.
\end{array}
\end{eqnarray}

The fifth $r$-matrix $r_5$ is the $r$-matrix of the Lorentz algebra, (\ref{L6}), and the
corresponding twist is given by the formula (\ref{L13}).

\end{document}